\DocumentMetadata{%
	pdfstandard=A-3u,
	pdfversion=1.7,
	lang=en-US,
}

\documentclass[balance,colorlinks,upint,subscriptcorrection,varvw,mathalfa=cal=boondoxo,spanish,german,vietnamese,russian,greek]{asmeconf}

\hypersetup{%
	pdfauthor={John H. Lienhard},									  
	pdftitle={ASME Conference Paper LaTeX Template},                  
	pdfkeywords={ASME conference paper, LaTeX template, BibTeX style},
	pdfsubject = {Describes the asmeconf LaTeX template},			  
	pdflicenseurl={https://ctan.org/pkg/asmeconf},
}

\usepackage{placeins}
\usepackage{color} 
\usepackage{soul}

\newcommand{\tl}[1]{\textcolor[rgb]{0,0.6,0.60}{#1}}
\usepackage[textsize=tiny,textwidth=0.5in]{todonotes}

\begin{document}


\ConfName{Proceedings of the ASME 2025 International Design Engineering Technical Conferences  \& \linebreak Computers and Information in Engineering Conference}
\ConfAcronym{IDETC/CIE2025}
\ConfDate{August 17-20, 2025} 
\ConfCity{Anaheim, California} 
\PaperNo{DETC2025-168547}


\title{Hybrid Powertrain Optimization for Regional Aircraft Integrating Hydrogen Fuel Cells and Aluminum Air Batteries} 

\SetAuthors{
	Harshal Kaushik\affil{1}, 
        Ali Mahboub Rad\affil{1}, 
        Korebami Adebajo\affil{2},
        Sobhan Badakhshan\affil{1}, 
	Nathaniel Cooper\affil{2}, \\
	Austin Downey\affil{2},  
	Jie Zhang\affil{1}\CorrespondingAuthor{jiezhang@utdallas.edu}
}

\SetAffiliation{1}{University of Texas at Dallas, Richardson, TX}
\SetAffiliation{2}{University of South Carolina, Columbia, SC}
 
%
%
%






\maketitle





\begin{abstract}

With the increasing demand for air travel and the urgency to reduce emissions, transitioning from fossil fuel-based propulsion systems is a critical step toward sustainable aviation. While batteries are widely used in urban air mobility, their long charging durations limit their feasibility for consecutive flights. Hybrid propulsion systems, which integrate fuel cells and batteries, offer a promising alternative due to their higher energy density and improved efficiency. This paper presents a novel hybrid powertrain architecture for regional aircraft, incorporating a hydrogen fuel cell, a lithium-ion battery, and an auxiliary aluminum-air battery. The proposed system is evaluated using real-world power demand data from a Cessna 208 aircraft. The hydrogen fuel cell acts as the primary power source, ensuring continuous operation, while the lithium-ion battery manages transient power fluctuations to enhance system stability. The aluminum-air battery is introduced as a high-energy emergency backup, providing extended endurance during critical situations. A mixed-integer optimization model is formulated for system sizing and power scheduling, ensuring optimal energy distribution among the power sources. Multiple operational scenarios are analyzed to evaluate system performance, particularly under emergency conditions, where power reliability is crucial. The results highlight the feasibility and effectiveness of the proposed hybrid architecture in improving energy efficiency and flight safety for regional aircraft applications.
\end{abstract}


\section{Introduction}

The shift toward sustainable, dependable energy sources is paramount in tackling global environmental challenges and mitigating climate change. The transportation sector significantly contributes to greenhouse gas emissions, accounting for roughly 2.5\% of global emissions; within this sector, aviation alone is responsible for nearly 12\% of transportation-related emissions, as reported in \cite{sparano2023hydrogen}. As air travel demand continues to rise, adopting cleaner and more efficient energy alternatives is vital to minimizing aviation’s environmental footprint and ensuring long-term sustainability.

Hydrogen has emerged as a promising alternative to fossil fuels in aviation, primarily due to its high gravimetric energy density. Research into hydrogen-based aviation has been ongoing for more than five decades, highlighting its potential as a cleaner energy carrier \cite{kadyk2018fuelcell}. Among hydrogen technologies, fuel cells stand out for their high efficiency, minimal noise, and zero harmful emissions, producing only water vapor \cite{A.V.Geliev2019}. However, as noted by Jarry et al. \cite{T.Jarry2021}, fuel cell systems face notable challenges such as slower dynamic response, limited power density, and accelerated aging, which hinder their standalone implementation in aviation.

Lithium-ion (Li-ion) batteries offer rapid response times, enhanced operational flexibility under dynamic conditions, and higher power density. However, their lower gravimetric energy density can compromise weight efficiency and reduce an aircraft’s payload capacity \cite{li2021optimal}. Numerous studies have examined all-electric battery-powered aircraft, with most focusing on short-range commuter flights due to battery energy density limitations. The study conducted by \cite{anker2025feasibility} highlights power constraints as a major challenge, emphasizing payload reductions for short-term viability and the need for battery advancements for long-term feasibility. Similarly, Bærheim et al. \cite{Baerheim2023battery} explore battery-electric aviation for regional flights up to 400~km, finding that while new aircraft designs meet operational needs, retrofitted models require improved battery energy density. In addition, Ebersberger et al. \cite{EBERSBERGER} assess an all-electric propulsion system for small commuter aircraft, analyzing power distribution, battery technologies, and energy storage while addressing efficiency, reliability, and safety considerations.

Fuel cell-battery hybridization, integrating both fuel cell and battery technologies into in a single power system has emerged as a promising approach to leverage their strengths while mitigating their weaknesses \cite{donateo2024retrofitting}.However, hybridization introduces challenges in energy management and sizing methodologies. As discussed in \cite{donateo2022conceptual}, these challenges can be formulated as an optimization problem, with sizing strategies and predictive algorithms ensuring seamless power allocation under varying operational conditions. 

Energy management and sizing optimization of fuel cell-battery hybrid-powered aircraft have received extensive research. For instance, Donateo et al. \cite{donateo2024retrofitting} investigate the retrofitting of ultralight aircraft with a fuel cell, lithium battery, and hydrogen storage system, optimizing hybridization ratios and comparing charge-depleting and charge-sustaining modes to assess their impact on fuel cell performance and hydrogen consumption. A refined sizing approach for fuel cell-battery hybrid systems is presented in \cite{park2022refined}, introducing a weight estimation model and conducting sensitivity analyses to determine the optimal power balance for enhanced aircraft performance. The study by Lei et al. \cite{lei2019energy} classifies energy management strategies into rule-based, intelligent-based, and optimization-based approaches, evaluating their efficiency and feasibility while analyzing trade-offs between fuel cells, batteries, and ultra-capacitors for unmanned aerial vehicle applications. Later in \cite{massaro2024optimal}, a comprehensive mathematical model has been formulated to optimize fuel cell performance across different scales, highlighting the importance of balancing hydrogen storage and electrochemical efficiency to enhance system reliability. Additionally, Hoenicke et al. \cite{hoenicke2021power} introduce a novel Power Management module to optimize energy distribution in hybrid propulsion, ensuring efficient fuel cell operation during cruise and battery support during high-power phases such as takeoff. Furthermore, this module integrates passive and active charging strategies to reduce battery size, improve reliability, and enhance overall safety.

Emergency energy reserves remain a challenge for hybrid fuel cell–battery systems because Li-ion batteries add significant weight when used solely for reserve capacity, and hydrogen-based systems do not scale effectively for this purpose. Consequently, increasing reserve capacity often encroaches upon the already limited usable range of these aircraft. In addition to hydrogen and battery-based hybrid solutions, alternative energy resources are now being explored for aviation and broader transportation applications. Among these, Aluminum–air (Al-air) batteries have garnered attention due to their high energy density, lightweight design, and long duration power supply capabilities. Their ability to deliver reliable emergency power makes them especially suitable for main energy source failures or divergence missions, where extended flight time or safe landing requires additional power. As noted in \cite{yang2002design}, Al-air batteries offer significant weight savings and a higher theoretical energy density compared to Li-ion batteries, rendering them an appealing option for auxiliary aviation use. Furthermore, aluminum is approximately 4050 times more abundant than lithium in the Earth’s crust \cite{xue2025air}, suggesting a more sustainable and widely available resource for energy storage. Nonetheless, issues such as anode corrosion, hydrogen evolution side reactions, and electrolyte degradation can compromise efficiency and lifespan. Moreover, as pointed out in \cite{chen2024identification}, the non-rechargeable nature of Al-air batteries necessitates anode and electrolyte replacement, limiting their practicality for continuous operation in primary energy applications. Despite these drawbacks, ongoing advancements in battery management, material engineering, and electrolyte optimization have the potential to bolster their viability for backup power and emergency energy support in aviation \cite{Pawlak2025}. However, only a few studies have investigated energy management and sizing optimization for Al–air–powered aircraft. One such study \cite{Vegh2025} examines the design and optimization of short-range, Al–air–powered aircraft for regional transport, analyzing weight, operating costs, and performance trade-offs. Recognizing that Al–air batteries offer high specific energy but low power output, the authors propose hybridizing them with Li-ion batteries to overcome power constraints and improve operational feasibility.


This paper extends the optimization model proposed in \cite{massaro2024optimal, hoenicke2021power} by incorporating aluminum air battery as an emergency backup power source and optimizing for sizing and energy management. We investigate this integration and assess aircraft performance in critical scenarios.  Our key contributions are outlined as follows:
\begin{itemize}
    \item We formulate a mixed-integer programming (MIP) model for optimal system sizing and power management, ensuring efficient energy distribution and hybridization ratio optimization.
    \item We introduce the integration of an aluminum-air battery within a hybrid powertrain consisting of a hydrogen fuel cell and a Li-ion battery. The aluminum-air battery serves as a backup power source, ensuring continued operation in the event of a fuel cell failure or providing additional power during unexpected scenarios such as rerouting or operational challenges.
    \item We derive the power curves for the Cessna 208 based on real-time flight data provided by pilots. Utilizing these power profiles, we design three distinct scenarios: normal operation, reduced hydrogen availability, and flight diversion. Each scenario is analyzed to ensure the aircraft can safely complete its mission and land securely.
\end{itemize}

The remainder of the paper is structured as follows. Section \ref{sec:2} presents the MIP formulation for system sizing, optimal power management, and hybridization ratio. Section \ref{sec:3} gives the experimental set up and Section \ref{sec:4} explores numerical case studies, detailing different experiments and analyzing the corresponding results. Finally, the paper concludes with a discussion of key findings in Section \ref{sec:5}.



\section{Mathematical Model}\label{sec:2}
In this section, we explain the proposed mixed-integer programming model for optimal sizing and power management.

\subsection{Objective Function}
The objective of this study is to optimize the size of the hybrid powertrain to minimize costs while ensuring that the aircraft's weight constraints are not exceeded. The analysis is based on the weight and volume constraints of the Cessna 208, which will be discussed in detail later in Section 3. Considering the removal of the original Pratt and Whitney PT6A-114A turboprop engine and the existing fuel storage, the available weight allocation for the fuel cell and battery system is estimated at 1200 kg. This serves as a key design constraint, defining the maximum allowable weight for integrating the hybrid power system.

Next, we define the design variables, which represent the energy requirements for the flight. The objective is to determine the optimal sizing of the hydrogen tank (L), fuel cell capacity (kWh), Li-ion battery capacity (kWh), and aluminum-air battery capacity (kWh). These variables are denoted as $V_{\text{H}}$, $E_{\text{fc}}$, $E_{\text{Li}}$, and $E_{\text{Al}}$, respectively. The first component of  objective function  is formulated as the minimization of the total weight of all hydrogen storage and power generation components. Thus, the optimization problem aims to minimize the following:
\begin{align} 
C_{\text{H}}^{\text{wt}} \ V_{\text{H}} + C_{\text{FC}}^{\text{wt}} \ E_{\text{fc}} + C_{\text{Li}}^{\text{wt}} \ E_{\text{Li}} + C_{\text{Al}}^{\text{wt}}\  E_{\text{Al}}, 
\end{align}
where the weight coefficients are finalized for batteries and fuel cell and listed in TABLE \ref{tab:weight_coeff}. 
\begin{table}[h]
\centering
\caption{Weight coefficient values}
\label{tab:weight_coeff}
\begin{tabular}{|l|l|l|l|l|}
\hline
\textbf{Coefficient} & \begin{tabular}[c]{@{}l@{}}$C_{\text{H}}^{\text{wt}}$ \\ (kg/L)\end{tabular} & \begin{tabular}[c]{@{}l@{}}$C_{\text{fc}}^{\text{wt}}$ \\ (kg/kWh)\end{tabular} & \begin{tabular}[c]{@{}l@{}}$C_{\text{Li}}^{\text{wt}}$\\ (kg/kWh)\end{tabular} & \begin{tabular}[c]{@{}l@{}}$C_{\text{Al}}^{\text{wt}}$\\ (kg/kWh)\end{tabular} \\ \hline
\textbf{Value}    & 1/11000   & 1.5   & 4   & {0.1234} \\ \hline
\end{tabular}
\end{table}

A second component of objective function is designed to collectively minimize the power expenditure of both the fuel cell and the batteries at every time instance $t$. The goal is to ensure that power consumption is minimized while meeting all operational constraints, thus maximizing the system's overall efficiency. Following is the second part of the objective function, where we seek to minimize the total power expenditure for both the fuel cell and the batteries across all time instances, formulated as:
\begin{align}
    \sum_{t \in T} P_{\text{fc}}^t + P_{\text{Li}}^t + P_{\text{Al}}^t.
\end{align}

{While the power demand profile is indeed a fixed input derived from real flight data and must be fully satisfied at every time step (enforced by the power balance constraint Equation \eqref{eq:power_demand_satify}), we include the power supplied by the fuel cell, lithium-ion battery, and aluminum-air battery in the objective function to enable optimal energy management. This formulation allows the model to strategically allocate power among the available sources, selecting combinations that not only satisfy demand but also minimize fuel consumption and enhance overall system efficiency. This integrated objective supports the design of a hybrid powertrain that is both energy-efficient and operationally sustainable.}

\subsection{Weight Constraints}
Next, we introduce the constraint sets, beginning with the weight constraint. This constraint ensures that the total weight of the combined power generation components and the hydrogen tank remains within the allowable weight limit. The hybrid energy system in the Cessna 208 is set to have the maximum allowable weight of 1200 kg.
\begin{align}
    C_{\text{H}}^{\text{wt}} \ V_{\text{H}} + C_{\text{FC}}^{\text{wt}} \ E_{\text{fc}} + C_{\text{Li}}^{\text{wt}} \ E_{\text{Li}} + C_{\text{Al}}^{\text{wt}}\  E_{\text{Al}} \leq 1200,
\end{align}
where the weight coefficients are mentioned earlier in TABLE \ref{tab:weight_coeff}. 
{As we described earlier, even though the objective function aims to minimize the total weight of the proposed hybrid powertrain, we explicitly include this maximum weight constraint to ensure compatibility with the baseline configuration of the Cessna 208, which is originally powered by a Pratt \& Whitney PT6A-114A turboprop engine. This constraint ensures that our replacement powertrain of hydrogen fuel cell and battery does not exceed the weight capacity allocated to the original propulsion system, thereby preserving the aircraft’s structural integrity, center-of-gravity limits, and aerodynamic balance. As the Cessna 208 airframe and performance characteristics are tightly coupled with its existing powertrain weight, our design prioritizes a modular retrofit without making significant modifications to the aircraft body. The weight constraint thus serves as a critical boundary condition, enforcing feasibility within the physical and regulatory limits of the original airframe.}


\subsection{Power Demand Satisfaction Constraint}
This constraint focuses on scheduling the power distribution to ensure that load requirements are met at every time instance while preventing any power generation unit, including the fuel cell and Li-ion battery, from being overloaded. The fuel cell and Li-ion battery serve as the primary hybrid power sources, whereas the aluminum-air battery functions as a reserved power supply, utilized only in emergency scenarios. The decision variables for this stage are $P_{\text{fc}}^\text{t}$, $P_{\text{Li}}^\text{t}$, and $P_{\text{Al}}^\text{t}$, representing the power generated by the fuel cell, Li-ion battery, and aluminum-air battery at each time instance $t$, respectively. The total flight duration is denoted as $T$, which, for simplicity, represents the set of all discretized time steps ranging from 0 to the end of the flight.

We begin by formulating the constraint that ensures the total power demand is met by the available power sources, namely the fuel cell, Li-ion battery, and aluminum-air battery. Mathematically, this condition can be expressed as follows:
\begin{align}\label{eq:power_demand_satify}
    P_{\text{fc}}^t + P_{\text{Li}}^t + P_{\text{Al}}^t = P_{\text{dem}}^t, \quad \forall t\ \in \ T,
\end{align}
where $P_{\text{dem}}^t$ is the demand at an instant $t$. The power demand data for the Cessna 208 has been gathered from the fuel flowdata provided by the pilots. By incorporating all these instances, we ensure that the power output sufficiently meets the load requirements under normal operating conditions. For exceptional extreme scenarios, a detailed discussion is provided later to assess and guarantee that the aircraft possesses adequate energy reserves for safe landings.

\subsection{Fuel Cell Constraints}
The third set of constraints ensures that the available hydrogen supply is adequate to sustain the fuel cell's energy requirements throughout the entire flight duration.
\begin{align}
    V_{\text{H}} \ \eta_{\text{fc}} \ H_{\text{LHV}} \ H_{\text{mass}} \geq E_{\text{fc}},
\end{align}
where $\eta_{fc}$ is the efficiency of fuel cell, $H_{\text{LHV}}$ is the lowest heating value of hydrogen, and $H_{\text{mass}}$ 
is the mass to volume ratio for hydrogen. Values for these are mentioned in TABLE \ref{tab:FC_properties}. 
\begin{table}[h]
\centering
\caption{Hydrogen fuel cell properties}
\label{tab:FC_properties}
\begin{tabular}{|l|l|l|l|}
\hline
\textbf{Coefficient} & \begin{tabular}[c]{@{}l@{}}$\eta_{\text{fc}}$ \\ \end{tabular} & \begin{tabular}[c]{@{}l@{}}$H_{\text{LHV}}$ \\ (kWh/kg)\end{tabular} & 
\begin{tabular}[c]{@{}l@{}}$H_{\text{mass}}$\\ (kg/L)\end{tabular}  \\ \hline
\textbf{Value}    & 0.55   & 33.33   & 0.09  \\ \hline
\end{tabular}
\end{table}

Next, we impose a constraint on the ramp rate of the fuel cell. We assume that there is a ramp rate limiting the change in the fuel cell's power output from $P_{\text{fc}}^{t_1}$ to $P_{\text{fc}}^{t_2}$. For the current scenario, we assume that within the time interval $\Delta t$, the fuel cell can increase its power by $0.1 E_{\text{fc}}$. This can be formulated as: 
\begin{align}
    P_{\text{fc}}^{t_1} - P_{\text{fc}}^{t_2} \leq 0.1 \ E_{\text{fc}}, \quad \forall t_1, t_2 \in T,
\end{align}
where $t_1$ and $t_2$ represent consecutive time instances within $T$. 
\subsection{Battery Constraint}
The state of charge (SOC) constraint on the Li-ion battery ensures that the battery's energy levels are properly tracked throughout the flight. This constraint is essential for maintaining energy balance and preventing over-discharge or overcharging, which could compromise the performance and longevity of the battery. Mathematically, it enforces that the energy stored in the battery at any time step $t_1$ is equal to the energy available at a subsequent time step $t_2$, adjusted for the power drawn from the battery over that interval.
\begin{align}
    SOC_{\text{Li}}^{t_1} \ E_{\text{Li}} = SOC_{\text{Li}}^{t_2} \ E_{\text{Li}} - P_{\text{Li}}^{t_1} \ \Delta t, \quad \forall t_1, t_2 \in T.
\end{align} 
Here, $P^{t_1}_{\text{Li}}$ represents the power supplied by the Li-ion battery at time $t_1$, while $\Delta t$ accounts for the duration over which this power is drawn. Similarly, the SOC equation for the aluminum-air battery is formulated in the same manner, establishing a relationship between $SOC^{t_1}_{\text{Al}}$, $E_{\text{Al}}$, and $P_{\text{Al}}^{t_1}$. This ensures that the aluminum-air battery's energy dynamics are accurately represented, enabling its role as an auxiliary power source during critical flight scenarios. 

\begin{figure} [bht!]
    \centering
    \includegraphics[width=0.99\linewidth]{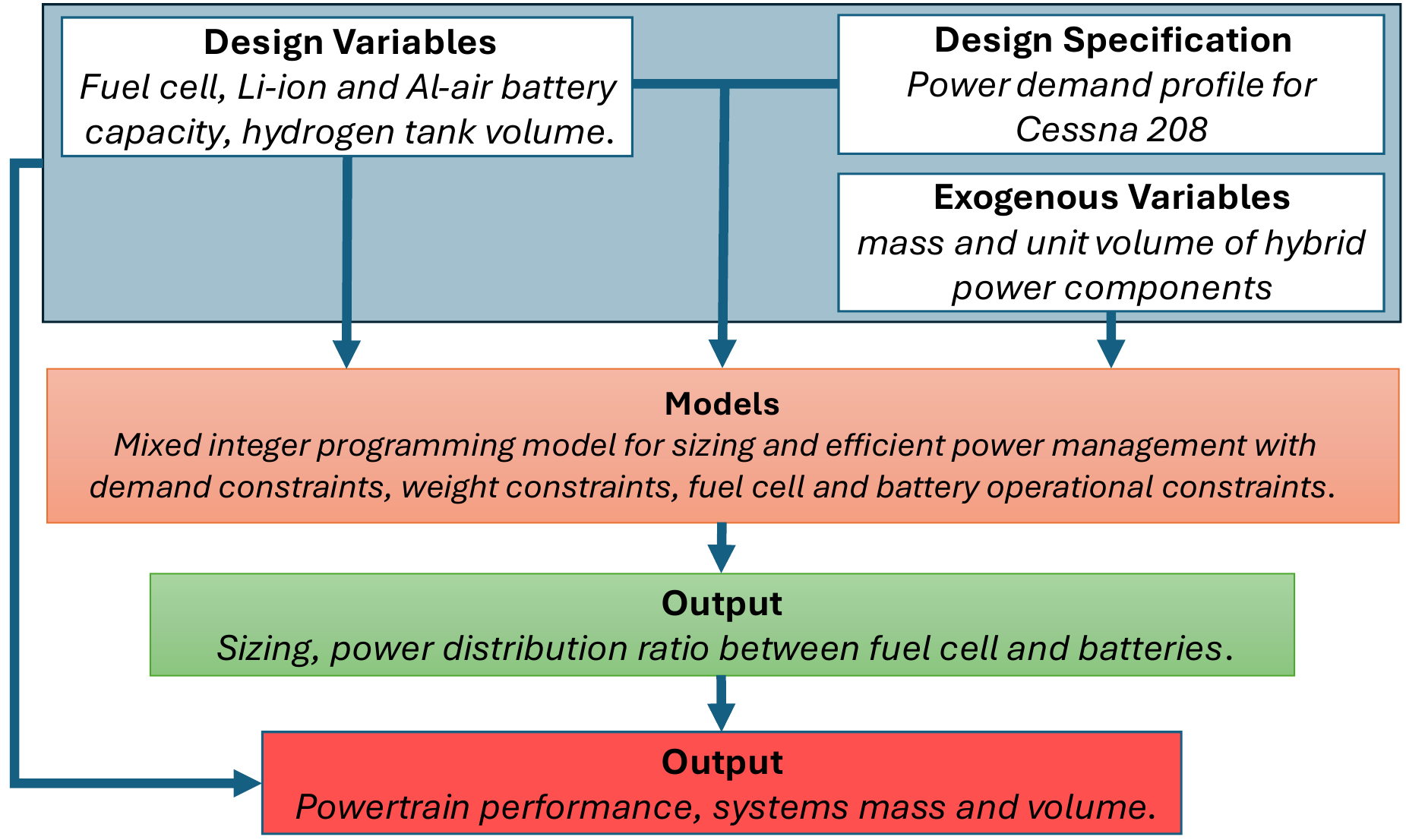}
    \caption{The proposed algorithm for optimal sizing and power management.} 
    

    \label{fig:flowchart_method}
\end{figure}

{We enforce the battery charge/discharge rate constraints implicitly by adopting a 1C rate limit for both the lithium-ion and aluminum-air batteries. Technically, a 1C rate implies that a battery can be fully charged or discharged within one hour. This constraint ensures that the battery’s power output or input at any time stays within its nominal capacity.}
The proposed methodology is compactly illustrated in the flowchart shown in Fig. \ref{fig:flowchart_method}.


\begin{figure} [t]
    \centering
    \includegraphics[width=\linewidth]{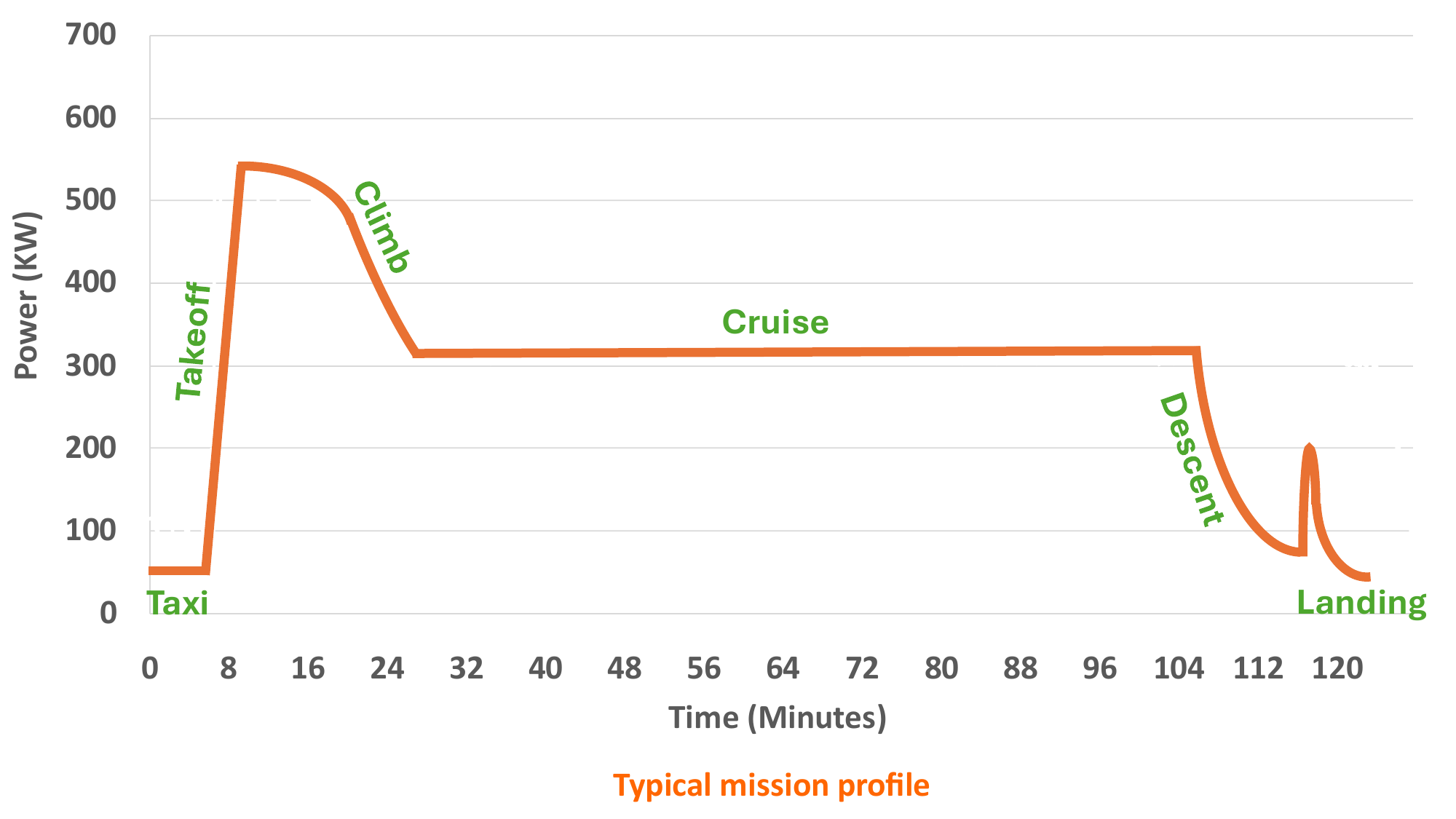}
    \includegraphics[width=\linewidth]{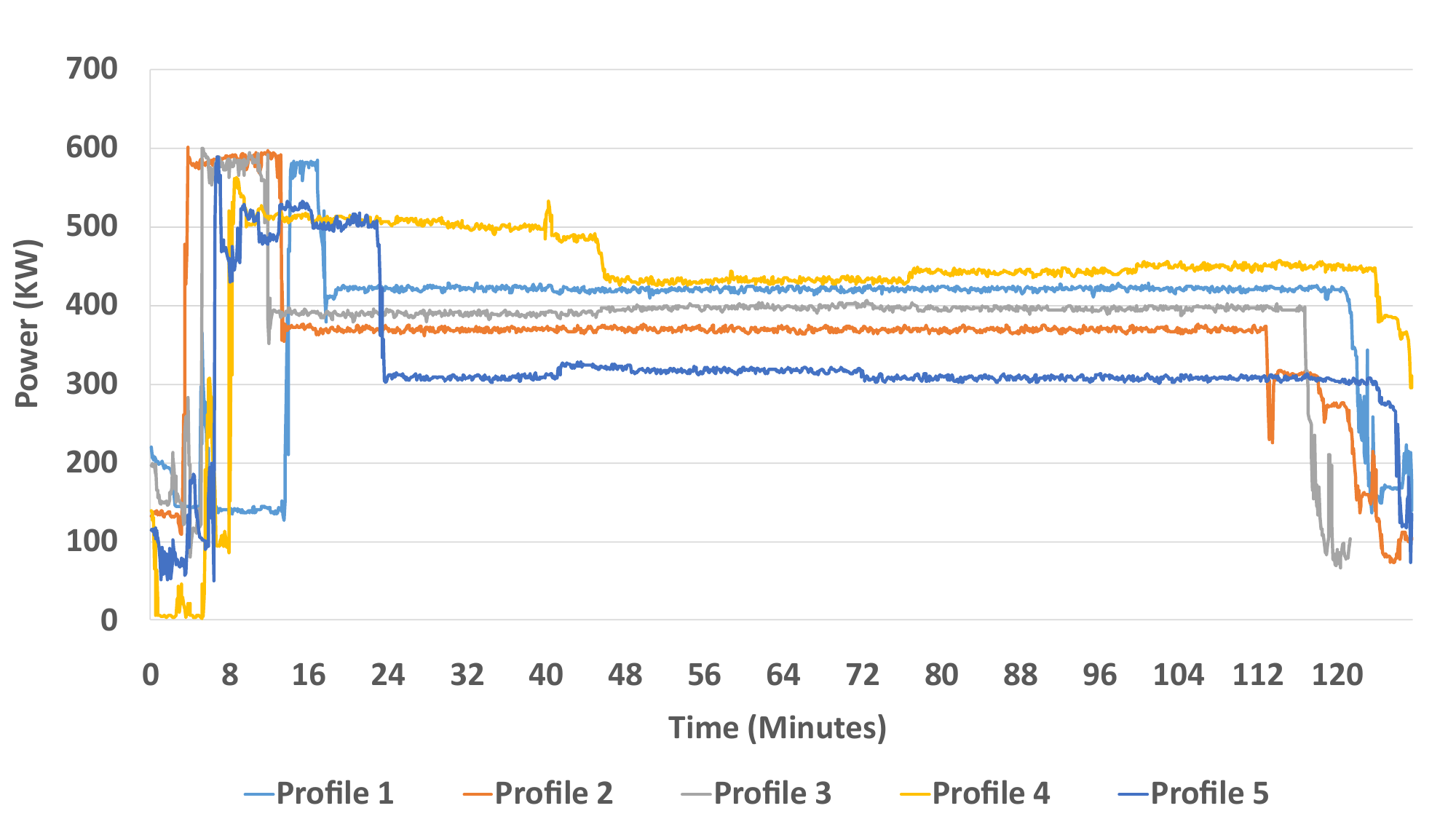}
    \caption{Cessna 208 power profiles from Columbia to Richmond. }  
    \label{fig:Columb_Rich_PowerProf_Compiled}
\end{figure}

\section{Experimental Set-up}\label{sec:3}
 In this section we introduce the selected aircraft (Cessna 208) followed by a detailed discussion of the power curves, including the methodology used to obtain them.


\subsection{Aircraft Specifications}
We utilize power curves derived from the Cessna 208 aircraft for a Columbia South Carolina, USA to Richmond  North Carolina, USA flight, which has an approximate duration of 130 minutes \cite{arts_lab_flight_data}. Our objective is to power the same aircraft with a zero-emission propulsion system using a hydrogen fuel cell and batteries while maintaining its aerodynamic performance, shape, size, and weight characteristics. To achieve this, we replace the conventional engine and fuel system with a hydrogen fuel cell and a combination of Li-ion and aluminum-air batteries, ensuring that the fundamental aircraft design remains unchanged.

The Cessna 208 has a fuel capacity of approximately 1000 kg and an engine (PT6A-114A) weighing around 200 kg. By eliminating these components, we allocate a total of 1200 kg to integrate the fuel cell, hydrogen storage, Li-ion battery, and aluminum-air battery required to sustain the entire flight. The fuel cell serves as the primary power source due to its high power density, efficiently driving the propeller. However, since fuel cells struggle with dynamic load variations, a Li-ion battery is incorporated to handle transient power demands. Additionally, an aluminum-air battery is included as an emergency energy source, ensuring sufficient power for a safe landing at the nearest airport in the event of a fuel cell failure. The precise capacities of all components will be determined through the experiments, explained next.


\subsection{Power Curves}
Multiple power curves have been obtained for the Columbia-to-Richmond flight, specifically for a conventionally powered Cessna 208 aircraft, illustrated in Fig. \ref{fig:Columb_Rich_PowerProf_Compiled}. These power curves capture various flight phases, including takeoff, climb, cruise, descent, and taxi. To facilitate optimized power system sizing and scheduling, we consider an average power curve that represents the overall power demand throughout the flight. 


Information on the power curves for various flights was obtained from fuel flow data recorded and provided by aircraft pilots. This data, sourced from Flightdata \cite{flightdata2024}, was converted to power used by the engine. The specific energy of the fuel commonly used in piston-powered aircraft was used to calculate the energy produced from the fuel consumed. By factoring in the engine’s efficiency at different points in time, this energy was then converted into usable power for the aircraft. During taxi, cruise, and descent, efficiency was assumed to be 30\%, while during takeoff, it was assumed to be 19\%, as these engines are typically less efficient at high power demands, such as during takeoff.

Due to the dynamic power variations during takeoff and climb, batteries play a crucial role in meeting the high transient power demands during these phases. Once the power requirement stabilizes during cruise, the fuel cell becomes the primary power source, while the battery usage is adjusted accordingly. The Li-ion battery is utilized to smooth out power fluctuations, ensuring stable operation, whereas the aluminum-air battery primarily serves as a backup energy source for safety, providing emergency power in the event of unexpected failures. 

In the following subsections, we conduct experiments to analyze the aircraft powertrain’s response under varying flight conditions and determine the required component sizing. Three experimental case studies are considered: the first examines a configuration utilizing only hydrogen and Li-ion batteries; the second evaluates powertrain performance under low hydrogen volume conditions; and the third investigates the aircraft’s ability to execute a safe landing while accounting for rerouting requirements.


\section{Numerical Experiments and Results}\label{sec:4}
\begin{figure} [t]
    \centering
    \includegraphics[width=\linewidth]{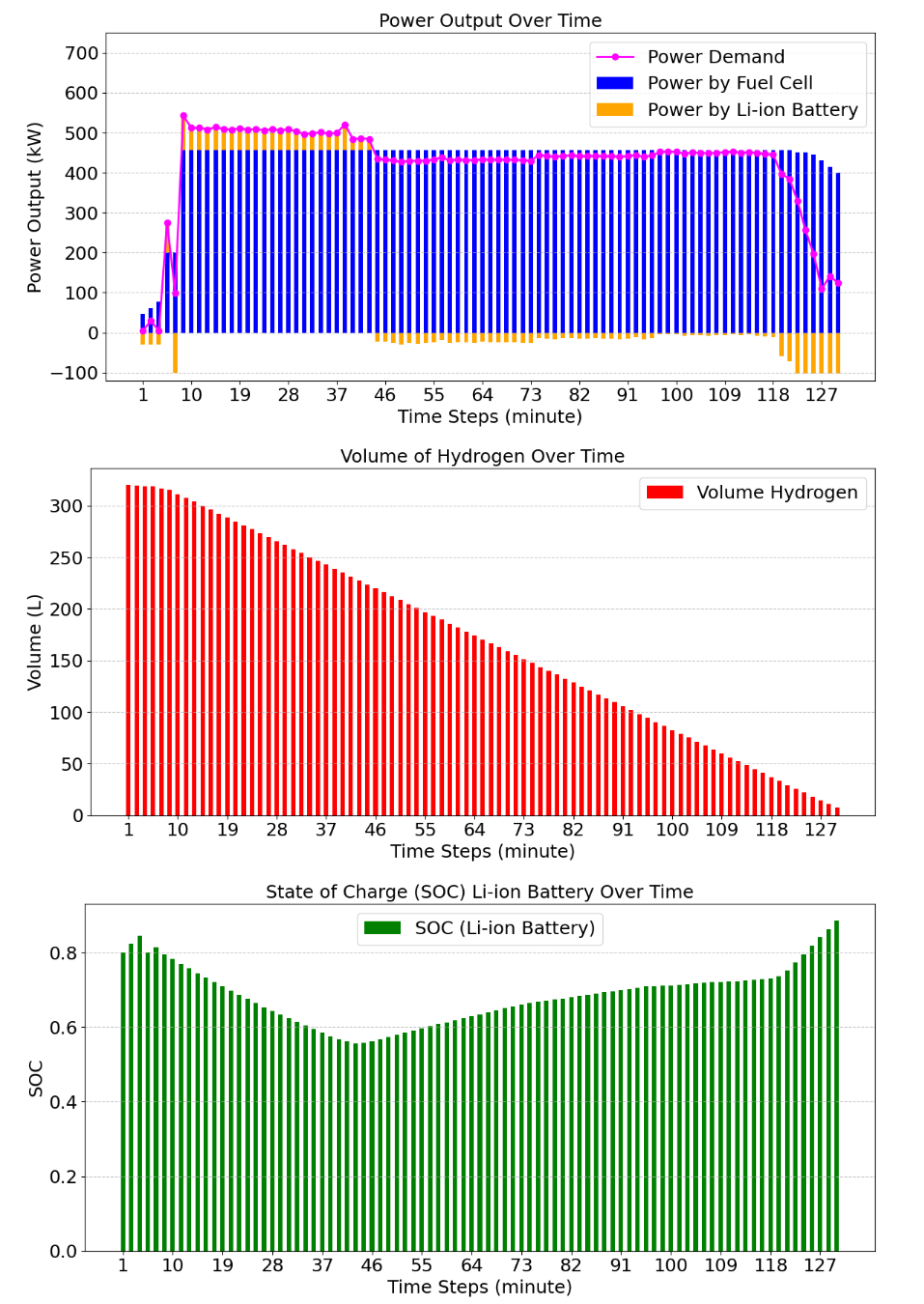}
    \caption{Simulation results of Experiment 1} 
    \label{fig:expt1}
\end{figure}

\begin{figure} [t]
    \centering
    \includegraphics[width=\linewidth]{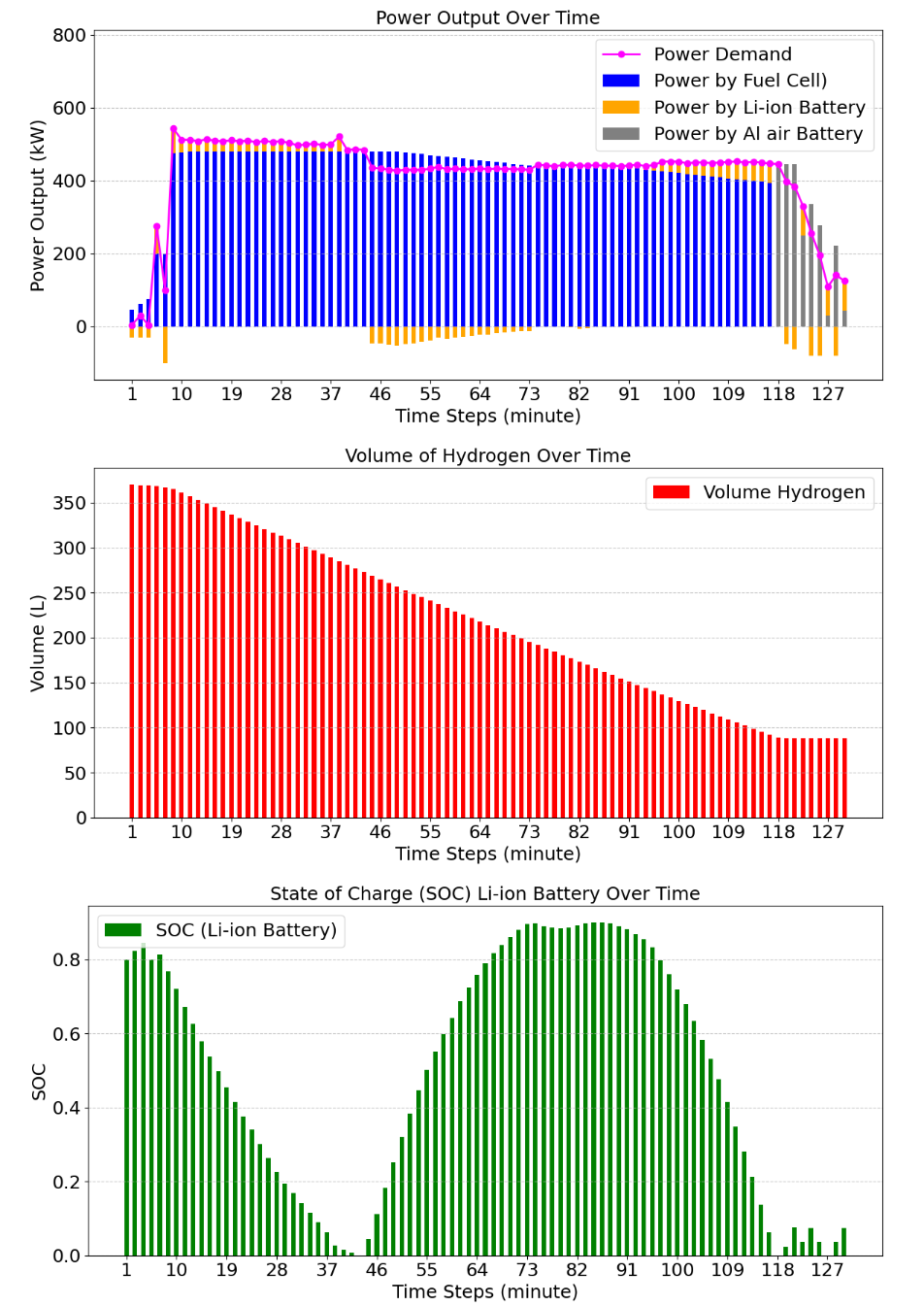}
    \caption{Simulation results of Experiment 2}
    \label{fig:expt2}
\end{figure}

\begin{figure} [t]
    \centering
    \includegraphics[width=\linewidth]{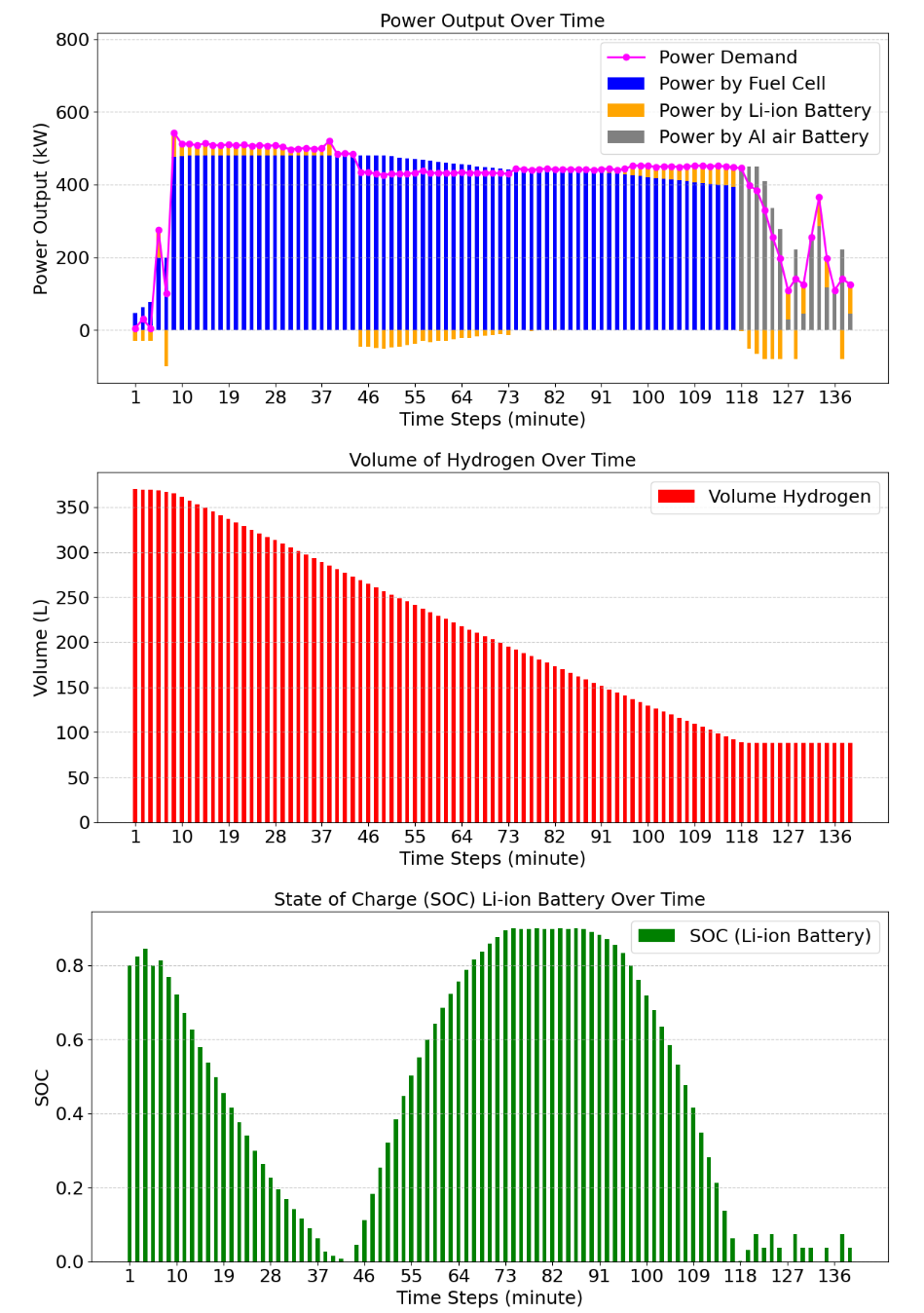}
    \caption{Simulation results of Experiment 3}
    \label{fig:expt3}
\end{figure}
In this section, we present the simulation experiments conducted using a real aircraft dataset as well as outline the experimental setup and discuss results.

\subsection{Experiment 1: Li-ion Battery Handling the Power Demand Fluctuations}
In this experiment, we use the average power curve derived from the power curves dataset shown in Fig. \ref{fig:Columb_Rich_PowerProf_Compiled}. The power is produced with a combination of a hydrogen fuel cell and a Li-ion battery. The fuel cell serves as the primary power source, supplying the majority of the power demand, particularly during the climb and cruise phases, where power requirements remain relatively stable. Minor perturbations in power demand are efficiently managed by the battery during those phases. Fig. \ref{fig:expt1}(a) illustrates the power output ratio between the fuel cell and the battery. During takeoff, there is a significant surge in power demand, which is primarily handled by the Li-ion battery. During the descent phase, power requirements are minimal and are efficiently met through the combined operation of the fuel cell and the battery. It is important to note that the fuel cell's deceleration response is relatively slow; therefore, any excess power generated is utilized to recharge the battery. Additionally, Fig. \ref{fig:expt1}(b) depicts the depletion of hydrogen storage over time, while Fig. \ref{fig:expt1}(c) presents the state of charge (SOC) of the battery throughout the flight. Table \ref{tab:exp_results_sizing} presents the sizing results for the fuel cell, hydrogen storage vessel, and Li-ion battery.

\begin{table*}[tbh!]
    \centering
        \caption{Sizing configurations for hybrid energy system components} \label{tab:exp_results_sizing}
    \renewcommand{\arraystretch}{1.2}
    \begin{tabular}{|l|c|c|c|c|c|}
        \hline
        \textbf{Experiment} & \textbf{Fuel Cell} & \textbf{Li-ion Battery} & \textbf{Aluminum-air} & \textbf{Hydrogen} & \textbf{Powertrain} \\
         & \textbf{capacity (kWh)} & \textbf{capacity (kWh)} & \textbf{battery capacity (kWh)} & \textbf{volume (L)} & \textbf{weight (kg)} \\
        \hline
        Experiment 1 & 458 & 100 & 0 & 320 & 1200 \\
        Experiment 2 & 480 & 80  & 450 & 370 & 1142 \\
        Experiment 3 & 480 & 80  & 450 & 370 & 1142 \\
        \hline
    \end{tabular}
\end{table*}

\subsection{Experiment 2: Aluminum-air Battery Activation at Low Hydrogen Levels}
In this experiment, we analyze a scenario in which the hydrogen volume reaches a critically low level, posing a potential risk to the flight. The key question addressed is how to ensure a safe flight continuation in such a situation. This is where the aluminum-air battery plays a crucial role. For this experiment, we assume that a safe hydrogen volume threshold is greater than 10\% of the total capacity (the aluminum-air battery is triggered below 100 L here in Figs. \ref{fig:expt2} and \ref{fig:expt3} for better illustration).

As in the previous case, power is initially supplied through a combination of a hydrogen fuel cell and a Li-ion battery. The fuel cell serves as the primary power source during the climb and cruise phases, where power demands remain relatively stable, while minor fluctuations are managed by the Li-ion battery. 

However, once the hydrogen volume drops below 10\% of the tank capacity, the aluminum-air battery is activated, and the fuel cell is shut down. At this stage, the aluminum-air battery becomes the primary power source, ensuring flight continuation and a safe landing. Furthermore, Fig. \ref{fig:expt2}(b) illustrates the depletion of hydrogen storage over time, while Fig. \ref{fig:expt2}(c) presents the SOC of the Li-ion battery throughout the flight. Table \ref{tab:exp_results_sizing} presents the sizing results for the fuel cell, hydrogen storage vessel, Li-ion battery, and aluminum-air battery.

\subsection{Experiment 3: Flight Rerouting Analysis}
In this third experiment, our primary objective is to analyze the flight behavior in the event of a rerouting scenario. There are situations where the aircraft may need to execute a go-around maneuver and return for landing due to missed approach, low visibility, or adverse environmental conditions. 

As in the previous cases, the takeoff, climb, and cruise phases are powered by the combined operation of the fuel cell and the Li-ion battery, as shown in Fig. \ref{fig:expt3}(a). However, during rerouting, the aircraft descends with lower hydrogen reserves, and the aluminum-air battery is activated to meet the additional energy demands required for the diversion. The aluminum-air battery is capable of sustaining this extra operation efficiently. Furthermore, Fig. \ref{fig:expt3}(b) illustrates the depletion of hydrogen storage over time, while Fig. \ref{fig:expt3}(c) presents the SOC of the Li-ion battery throughout the flight.

\subsection{Optimal Sizing Results}
The sizing outcomes for the fuel cell, hydrogen storage vessel, Li-ion battery, and aluminum-air battery are summarized in Table \ref{tab:exp_results_sizing}. These results represent the optimum quantity of lithium ion battery energy, fuel cell energy, hydrogen volume, and aluminium air battery energy required to complete the flight in the emergency situations. These quantities are selected by obeying the weight constraints, sizing constraints, and operational constraints.

\FloatBarrier

\section{Conclusion}\label{sec:5}
A comprehensive analysis of a hybrid powertrain for a regional aircraft, i.e., the Cessna 208, integrating a hydrogen fuel cell, a Li-ion battery, and an aluminum-air battery is performed to enhance operational resilience. A mixed-integer programming model is formulated for optimal system sizing and power management, ensuring efficient energy distribution and hybridization optimization. Through the use of real-world flight data, we evaluate power demands under various scenarios, including nominal operation, reduced hydrogen availability, and flight diversion.  Findings highlight the critical role of the aluminum-air battery in ensuring flight safety during emergency conditions, particularly in cases of fuel cell failure or unexpected rerouting.
This study provides valuable insights into powertrain sizing and scheduling, with plans to expand the dataset by incorporating additional flight profiles for enhanced robustness. 

Potential future work will involve developing a stochastic sizing model incorporating diverse origin-destination datasets to enhance hybrid-electric propulsion for regional aircraft.


\section*{Acknowledgments}
This material is based in part upon work supported by the Air Force Office of Scientific Research (AFOSR) through award no. FA9550-21-1-0083. This work is also partly supported by the National Science Foundation (NSF) grant number 2237696. This project is also partially supported by the University of South Carolina Office of Undergraduate Research through the Magellan Scholar Program. Any opinions, findings, conclusions, or recommendations expressed in this material are those of the authors and do not necessarily reflect the views of the National Science Foundation, the United States Air Force, or the University of South Carolina.



\bibliographystyle{IEEEtran}
\bibliography{references}

\end{document}